\documentclass{amsart}
\usepackage[latin1]{inputenc}
\usepackage[english]{babel}
\usepackage{amssymb,amsmath,amsthm}
\usepackage{amsfonts}
\usepackage{url}
\usepackage{enumerate}
\theoremstyle{plain}
\newtheorem{theorem}{Theorem}
\newtheorem{corollary}{Corollary}
\newtheorem{remark}{Remark}
\newtheorem{problem}{Problem}

\newtheorem{lemma}{Lemma}
\newtheorem{proposition}{Proposition}
\newtheorem{fact}{Fact}
\newtheorem{claim}{Claim}

\newcommand{\N}{\mathbb{N}}

\newcommand{\B}{\mathbb{B}}

\begin{document}

\title{On a quasi-ordering  on Boolean functions}

\author{Miguel Couceiro}
\address{Department of Mathematics, Statistics and Philosophy\\
University of Tampere\\
Kalevantie 4, 33014 Tampere, Finland}
\email{Miguel.Couceiro@uta.fi}
\thanks{The work of the first author was partially supported by the Graduate School in Mathematical
 Logic MALJA, and by grant \#28139 from the Academy of Finland.
}
\author{Maurice Pouzet}
\address {PCS, Universit\'e Claude-Bernard Lyon1,
 Domaine de Gerland -b\^at. Recherche [B], 50 avenue Tony-Garnier, F$69365$ Lyon cedex 07, France}
\email{pouzet@univ-lyon1.fr }
\thanks {The work of the second named author was supported by INTAS}


\date{January, 2006}
\maketitle

\begin{abstract}  It was proved few years ago that classes of Boolean functions
 definable by means of functional equations \cite{EFHH}, or equivalently, by means of
relational constraints \cite{Pi2},
 coincide with  initial
 segments of the quasi-ordered set $(\Omega, \leq )$ made of the set $\Omega$ of Boolean functions,
 suitably quasi-ordered.
 The resulting ordered set $(\Omega/\equiv, \sqsubseteq)$ embeds into 
$([\omega]^{<\omega}, \subseteq)$,
 the set -ordered by inclusion- of finite subsets of the set $\omega$ of integers. We
 prove that $(\Omega/\equiv, \sqsubseteq)$ also embeds $([\omega]^{<\omega}, \subseteq)$.
 We prove that initial segments of $(\Omega, \leq )$ which are definable by finitely many 
obstructions coincide with classes defined by finitely many equations. This gives, 
in particular, that the classes of Boolean functions with a bounded number of essential variables
 are finitely definable. As an example, we provide a concrete characterization of the subclasses
 made of linear functions.
\end{abstract}

\section{Introduction}
Two approaches of Boolean definability have been considered recently. One in terms of
 functional equations \cite{EFHH}, an other in terms of relational constraints \cite{Pi2}.
 It turns out that these two approaches define the same classes of Boolean functions.
These classes have been completely described by means of a quasi-order on the set $\Omega$
 of all Boolean functions. The quasi-order is the following: for two functions $f,g\in \Omega$ set
$g\leq  f $ if $g$ can be obtained from $f$ by identifying, permuting or adding variables.
 These classes coincide with initial segments for this quasi-ordering
called \emph{identification minor} in \cite{EFHH}, \emph{minor} in \cite {Pi2}, \emph{subfunction}
in \cite {Z}, and \emph{simple variable substitution} in \cite{CF1}. 
 Since then, greater emphasis on this quasi-ordering has emerged. For an example, it was observed that    
  $\Omega$ is the union of four blocks with no comparabilities in between,
 each block made of the elements above a minimal element. 
 In \cite{Pi2}, Pippenger showed that $\Omega$ contains infinite antichains.
A complete classification of pairs $C_1, C_2$ of particular initial segments 
("clones") for which $C_2\setminus C_1$ contains no infinite antichains was given in \cite{C3}. 
Our paper is a contribution to the understanding  of this quasi-ordering. 

Some properties are easier to express in terms of the poset $(\Omega/\equiv, \sqsubseteq)$ associated with the
 quasi-ordered set $(\Omega, \leq )$ and  made of the equivalence classes associated with the equivalence 
$\equiv$ defined by $f\equiv g$ if $f\leq  g$ and $g\leq  f$. As we will see (Corollary  \ref {levels}),
 for each $x\in \Omega/\equiv$, the initial segment $\downarrow x:=\{y\in \Omega/\equiv: y\leq x\}$ 
is finite, hence $(\Omega/\equiv, \sqsubseteq)$ decomposes into the levels 
$\Omega/\equiv_{0}, \dots  \Omega/\equiv_{n}, \dots$, where $\Omega/\equiv_{n}$ is the set of minimal
 elements of $\Omega/\equiv\setminus \cup \{ \Omega/\equiv_{m}: m<n\}$. Moreover, each level is finite;
 for an example $\Omega/\equiv_{0}$ is made of four elements (the equivalence classes of the two constants 
functions, of the identity and of the negation of the identity). This fact leads to the following:
 \begin{problem}How does the map $\varphi_{\Omega/\equiv}$, which counts for every $n$ the number 
$\varphi_{ \Omega/\equiv}(n)$ of elements of $\Omega/\equiv_{n}$, behave?
 \end{problem}
 
From the fact that  for each  $x\in \Omega/\equiv$, the initial segment $\downarrow x$ is finite
it follows that initial segments  of $(\Omega/\equiv, \sqsubseteq)$ correspond bijectively to antichains of  $(\Omega/\equiv, \sqsubseteq)$. Indeed, for each antichain $A\subseteq (\Omega/\equiv,  \sqsubseteq)$, the set $Forbid (A):= \{y\in \Omega/\equiv \; : x\in A \Rightarrow x\not \sqsubseteq y\}$ is an initial segment of $(\Omega/\equiv, \sqsubseteq)$. Conversely, each  initial segment $I$ of  $(\Omega/\equiv, \sqsubseteq)$ is of this form  (if $A$ is the set of minimal elements of  $\Omega/\equiv \setminus I$, then since for each $x\in \Omega/\equiv$ the set $\downarrow x$ is finite, $I=Forbid (A)$).  Viewing the elements of $A$ as obstructions, this amounts to  say that {\it every initial segment can be defined by a minimal set of obstructions}.  

Another feature of this poset, similar in importance, is the fact that it is {\it up-closed},
 that is for every pair $x,y\in (\Omega/\equiv)$, the final segment $\uparrow x\cap \uparrow y$
 is a finite union (possibly empty) of final segments of the form  $\uparrow z$. This means that the
 collection of initial segments of the form $Forbid(A)$ where $A$ runs throught the finite antichains of 
 $\Omega/\equiv$ which is closed under finite intersections is also closed under finite unions. 

Such initial segments have a natural interpretation in terms of Boolean functions. Indeed, as we have said, initial segments of $(\Omega, \leq )$ coincide with  equational classes.  Each of these initial segments  identifies to an  initial segment of $(\Omega/\equiv, \sqsubseteq)$ and, as in this case, can be written  as  $Forbid (A)$ for some antichain $A$ of $(\Omega, \leq )$ (the difference with an initial segment of $(\Omega/\equiv, \sqsubseteq)$ is that the antichain $A$ is not unique). Let us   
 consider the set $ \mathcal F$ of classes which can be defined by finitely many equations.  
They are characterized by the following theorem.

\begin{theorem}\label{thm2}
For an initial segment $I$ of $(\Omega, \leq )$, the following properties are equivalent:
\begin{enumerate}[{(i)}]
\item $I\in \mathcal F$;
\item $I$ is definable by a single equation; 
\item $I= Forbid (A)$ for some finite antichain.
\end{enumerate}
\end{theorem}
The following lemma reassembles the main properties of $\mathcal F$. 
\begin{lemma}\label{mainlemma}
\noindent
\begin{enumerate}
\item $\mathcal F$ is closed under finite unions and finite intersections;
\item $Forbid(\{f\})\in \mathcal F$ for every $f\in \Omega$;
\item $\downarrow f\in \mathcal F$ for every $f\in \Omega$;
\item the class of $f\in \Omega$ with no more than $k$ essential variables belongs to $\mathcal F$ for every integer $k$. 
\end{enumerate}
\end{lemma}
The class of linear operations
(w.r.t the $2$-element field) belongs to $\mathcal F$; we give an explicit equation defining the class of linear operations with at most $k$ essential variables.
 Our proof makes use of basic linear algebra over the $2$-element field.

The set $\mathcal F$ ordered by inclusion is  a bounded distributive  lattice. 
As it is well known \cite{daveypriestley} a  bounded distributive lattice  $T$  is characterized by its
 {\it Priestley space}, that is the collection of prime filters of $T$, the {\it spectrum of T},
  ordered by inclusion and equipped with the topology induced by the product topology on $\mathfrak P(T)$.
 In our case, $\mathcal F$ is dually isomorphic to the sublattice of $\mathfrak P(\Omega/\equiv)$
 generated by the final segments of the form $\uparrow x$ for $x\in \Omega/\equiv$.
 This lattice is the {\it tail-lattice of $(\Omega/\equiv, \sqsubseteq)$}.
 From the fact that $(\Omega/\equiv, \sqsubseteq)$ is up-closed and has finitely many minimal elements, it 
follows that the Priestley space of the  tail-lattice of $(\Omega/\equiv, \sqsubseteq)$
 is the set $\mathcal J(\Omega/\equiv, \sqsubseteq)$ of ideals of $(\Omega/\equiv, \sqsubseteq)$ 
ordered by inclusion and equipped with the topology induced by the product topology on
 $\mathfrak P(\Omega/\equiv)$ (in \cite {BPZ}, Theorem 2.1 and Corollary 2.7). Hence we have:
\begin{theorem} The Priestley space of the lattice $\mathcal F$ ordered by reverse inclusion is the set
 $\mathcal J(\Omega/\equiv, \sqsubseteq)$ of ideals of $(\Omega/\equiv, \sqsubseteq)$ ordered by
 inclusion and equipped with the topology induced by the product topology on $\mathfrak P(\Omega/\equiv)$.
\end{theorem}

This result ask for a description of $\mathcal J(\Omega/\equiv, \sqsubseteq)$. 
We prove that it embeds the poset $(\mathfrak P(\omega), \subseteq )$, the power set of $\omega$,  ordered by inclusion. 

Our proof is a by-product of an attempt to locate $(\Omega/\equiv, \sqsubseteq)$ among posets, that we now describe.  
There are two well-known ways of
 classifying posets. One with respect to isomorphism, two posets $P$ and $Q$ being \emph{isomorphic}
 if there is some order-isomorphism from $P$ onto $Q$. The other w.r.t. equimorphism, $P$ and $Q$ being
 \emph{equimorphic} if $P$ is isomorphic to a subset of $Q$, and $Q$ is isomorphic to a subset of $P$.
 Given a  poset $P$, one may ask to which well-known poset $P$ is isomorphic or, if this is too difficult,
 to which $P$ is equimorphic. If $P$ is the poset $(\Omega/\equiv, \sqsubseteq)$, we cannot answer the
 first  question. We answer the second.

Let $[\omega]^{<\omega}$ be the set of finite subsets of the set $\omega$ of integers.
Once ordered by inclusion, this yields the poset $([\omega]^{<\omega}, \subseteq)$. This poset decomposes into levels, the  $n$-th level being made of the $n$-element subsets of $\omega$.  Since all its  levels (but one) are infinite, it is not isomorphic to $(\Omega/\equiv, \sqsubseteq)$.  But:
\begin{theorem}\label{main}  $(\Omega/\equiv, \sqsubseteq)$ is equimorphic to
$([\omega]^{<\omega}, \subseteq)$.
\end{theorem}

As it is well-known and easy to see, the poset $([\omega]^{<\omega}, \subseteq)$ contains an
isomorphic copy of every countable poset $P$ such that the initial segment $\downarrow x$
is finite for every $x\in P$. Since $(\Omega/\equiv, \sqsubseteq)$ enjoys  this property,
it embeds into $([\omega]^{<\omega}, \subseteq)$.
 The proof that $([\omega]^{<\omega}, \subseteq)$ embeds into $(\Omega/\equiv, \sqsubseteq)$ is based on a strenthening of a construction of an infinite antichain in $(\Omega, \leq )$ given in \cite {Pi2}. 

  Since $\mathcal J([\omega]^{<\omega} , \subseteq)$ is isomorphic to $(\mathfrak P(\omega), \subseteq)$, $\mathcal J(\Omega/\equiv, \sqsubseteq)$ embeds $(\mathfrak P(\omega), \subseteq)$, proving our claim above.  

 This work was done while the first named author visited the Probabilities-Combinatoric-Statistic  group at the Claude-Bernard University in Gerland during the fall of 2005. 

\section{Basic notions and basic results}
\subsection{Partially ordered sets and initial segments}

A \emph{quasi-ordered set} (qoset) is a pair $(Q,\leq )$ where $Q$ is an arbitrary set and
 $\leq $ is a \emph{quasi-order} on $Q$, that is, a reflexive and transitive binary relation on $Q$. If the quasi-order
is a \emph{partial-order}, i.e., if it is in addition antisymmetric, then this qoset is said to be a \emph{partially-ordered set} (poset). \emph{The equivalence $\equiv$  associated to  } $\leq$ is defined by $x\equiv y$ if $x\leq y$ and $y\leq x$. We denote $x<y$ the fact that $x\leq y$ and $y\not \leq x$.  We denote $\overline x$   the equivalence class of $x$ and $Q/\equiv$ the set of equivalence classes. The image of $\leq$ via the  quotient  map  from $Q$ into  $Q/\equiv$ (which associates $\overline x$ to $x$)  is an order, denoted $\sqsubseteq$. According to our notations, we have $x<y$ if and only if $\overline x  \sqsubset \overline y$. Throught this map,  properties of qosets translate into properties of posets.   The consideration  of a poset rather than a qoset is then matter of convenience.  

Let $(Q,\leq )$ be a qoset. A subset $I$ of $Q$ is  an \emph{initial segment} if it contains
every $q'\in Q$ whenever $q'\leq q$ for some $q\in I$.
We denote by $\downarrow X$
the initial segment \emph{generated by} $X\subseteq Q$, that is,
 \begin{displaymath}
\begin{array}{l}
\downarrow X=\{q'\in Q: q'\leq q\textrm{ for some }q\in X\}.
\end{array}
\end{displaymath}
If $X:=\{x\}$,  we use the notation $\downarrow x$ instead of $\downarrow \{x\}$. 
An initial segment of the form $\downarrow x$ is \emph{principal}. A \emph{final segment} of 
$(Q, \leq)$ is an initial segment for the dual quasi-order.
 We denote $\uparrow X$ the final segment generated by $X$ and use $\uparrow x$ if $X:= \{x\}$.
 Given a subset $X$ of $Q$, the set $Q\setminus \uparrow X$ is an initial segment of $Q$;
 we will rather denote it $Forbid (X)$ and refer to the members of $X$ as {\it obstructions}.
 We denote by $I(Q, \leq)$ the poset made of the initial segments of $(Q, \leq)$ ordered by inclusion. For an example $I(Q, =)=(\mathfrak P(Q), \subseteq)$. 
 An {\it ideal} of $Q$ is a non-empty initial segment $I$ of $Q$ which is {\it up-directed}, 
this condition meaning that for every  $x,y\in I$  there is some $z\in I$ such that $x,y\leq z$. 
We denote by $\mathcal J(Q, \leq)$ the poset made of the ideals of $(Q, \leq)$ ordered by inclusion. 

Let $(Q,\leq )$ and $(P,\leq )$ be two posets.
A map $e:Q\rightarrow P$ is an  \emph{embedding} of $(Q,\leq )$ into $(P,\leq )$
if satisfies the condition
 \begin{displaymath}
\begin{array}{l}
q'\leq q \textrm{ if and only if  }e(q')\leq e(q)
\end{array}
\end{displaymath}
Such a map is necessarily one-to-one. If it is surjective, this is an \emph{isomorphism} of $Q$ onto $P$. For an example $\mathcal J([\omega], \subseteq )$ is isomorphic to $(\mathfrak P(\omega), \subseteq)$.  
 
Hence an embedding of $Q$ into $P$ is an isomorphism of $Q$ onto its image. The relation $P$ \emph{is  embeddable into} $P$ if there is some embedding from $Q$ into $P$ is a quasi-order on the class of posets. Two posets  which are equivalent with respect to this quasi-order, that is which embed in each other are said \emph{equimorphic}.
We note that if $(Q, \leq)$ is a qoset the  quotient map from $Q$ onto $Q/\equiv$ induces an isomorphism from $I(Q, \leq)$ onto $I(Q/\equiv, \sqsubseteq)$ and from $\mathcal J(Q, \leq)$ onto $\mathcal J(Q/\equiv, \sqsubseteq)$. 

A \emph{chain}, or a \emph{linearly ordered set}, is
a poset in which all elements are pairwise comparable with respect to an order $\leq $.
 By an \emph{antichain} we simply mean a set of pairwise incomparable elements. 

Let $(P, \leq)$ be  a poset. Denote by  $Min (P)$ the subset of $P$ made of minimal elements of $P$.
 Define inductively the sequence $(P_n)_{n\in \N}$ setting $P_0:= Min (P)$ and 
$P_n:= Min( P\setminus \cup \{P_{n'}: n'<n\})$. For each integer $n$, the set $P_n$ is an antichain, called a {\it level} of $P$. If $P_n$ is non-empty, this is the 
\emph{$n$-th level} of $P$.  For $x\in P$,  we write $h(x,P)=n$ if $x\in P_n$.  Trivially, we have: 

\begin{lemma} $P$ is the union of the $P_n$'s whenever for every $x\in P$, the initial segment $\downarrow x$ is finite. 
\end{lemma}

We will need the following result. It belongs to the folklore of the theory of ordered sets.
 For sake of completeness we give a proof.

\begin{lemma} \label{finite} A poset $(P, \leq)$ embeds into $([\omega]^{<\omega}, \subseteq)$
 if and only if $P$ is countable and for every $x\in P$, the initial segment $\downarrow x$ is finite.
\end{lemma}
\begin{proof}The two conditions are trivially necessary.
 To prove that they suffice, set $\varphi (x):= \downarrow x$.
 This defines an embedding from $(P, \leq)$ into $([\omega]^{<\omega}, \subseteq)$.
\end{proof}

\subsection{Boolean functions}

Let $\mathbb{B}: = \{0, 1\}$. A \emph{Boolean function} is a map $f : \mathbb{B}^n \rightarrow \mathbb{B}$, for some positive integer $n$ called the \emph{arity} of $f$.  By a \emph{class} of Boolean functions,  we simply mean a set
 $K \subseteq \Omega $, where $\Omega $ denotes the set $\bigcup_{n \geq 1} \mathbb{B}^{\mathbb{B}^n}$
 of all Boolean functions. For $i,  n\in \N^*$ with
$i\leq n$, define the
$i$-{\it th $n$-ary projection}
$e^{n}_{i}$ by setting $e^{n}_{i}(a_{1},\dots,a_{n}):= a_{i}$. Set
$ I_c:= \{e^{n}_{i}: i,n \in \N^*\}$.
These {$n$-ary projection maps} are also called \emph{variables}, and denoted $x_1, \ldots, x_n$,
where the arity is clear from the context.
 If $f$ is an $n$-ary Boolean function and $g_1, \ldots, g_n$ are $m$-ary Boolean functions,
 then their \emph{composition} is the $m$-ary Boolean function $f(g_1, \ldots, g_n)$,
  whose value on every ${\bf a} \in \mathbb{B}^m$  is
$f(g_1({\bf a}), \ldots,\linebreak[0] g_n({\bf a}))$.
  This notion is naturally extended to classes $I,J\subseteq \Omega $, by defining
 their \emph{composition} $I\circ J$ as the set of all composites of functions in
$I$ with functions in $J$, i.e.
\begin{displaymath}
  I\circ J=\{f(g_1,\ldots ,g_n)\mid n,m\geq 1, f\textrm{ $n$-ary in $I$, }g_1,\ldots ,g_n
\textrm{ $m$-ary in $J$} \}.
\end{displaymath} When $I=\{f\}$,
we write $f\circ J$ instead of $\{f\}\circ J$.
Using this terminology, a \emph{clone} of Boolean functions is defined as a class $C$ containing all projections
and idempotent with respect to class composition, i.e., $C\circ C=C$.
As an example, the class $I_c$ made of  all projections is a clone.
For further extensions see e.g. \cite{CFL,CF1,CF2, CF3}.

\bigskip

An $m$-ary Boolean function $g$
is said to be obtained from an $n$-ary Boolean function $f$
 by \emph{simple variable substitution}, denoted $g\leq  f$,
if there are $m$-ary projections $p_1,\ldots ,p_n\in I_c$ such that
$g=f(p_1,\ldots ,p_n)$. In other words,
  \begin{displaymath}
g\leq  f \quad \text{if and only if }\quad  g\circ I_c\subseteq f\circ I_c.
\end{displaymath}
Thus $\leq  $ constitutes a quasi-order on $\Omega $. If $g\leq  f$ and $f\leq  g$, then $g$ and $f$ are said to be
\emph{equivalent}, $g\equiv f$. Let $\Omega /\equiv $ denote the set of all equivalent classes
of Boolean functions and let $\sqsubseteq $ denote the partial-order induced by $\leq  $.
A class $K\subseteq \Omega $
is said to be \emph{closed under simple variable substitutions}
if each function obtained from a function $f$ in ${K}$ by simple variable substitution
 is also in ${K}$.  In other words, the class ${K}$
 is closed under simple variable substitutions if and only if
$K/\equiv $ is an initial segment of $\Omega /\equiv $.
(For an early reference on the quasi-order $\leq  $ see e.g. \cite{W} and for futher background see \cite{EFHH, Pi2, Z, CF1, C2, C3}. 
For variants and generalizations see e.g. \cite{CF2,CF3,L,Le1,Le2}.)

\subsubsection{Essential variables and minors}
Let $f: \mathbb{B}^n\rightarrow \mathbb{B}$ be an $n$-ary Boolean function. For each $1\leq i\leq n$, $x_i$ is said to be an \emph{essential variable of $f$} if there are 
$a_1,\ldots ,a_{i-1},a_{i+1},\ldots ,a_n$ in $ \mathbb{B}$ such that 
\begin{displaymath} f(a_1,\ldots ,a_{i-1},0,a_{i+1},\ldots ,a_n)\not=f(a_1,\ldots ,a_{i-1},1,a_{i+1},\ldots ,a_n). 
\end{displaymath} 
Otherwise, $x_i$ is called \emph{a dummy variable of $f$}. The \emph{essential arity of} $f$, denoted $ess(f)$ is the number of its essential variables.   
Note that constant functions are the only Boolean functions whose variables are all dummy.  
\begin{lemma}\label{ess}
\noindent
\begin{enumerate}
\item If $g < f$  then $ess(g)< ess(f)$;
\item For every Boolean function $f$ we have
  \begin{displaymath}
max\{ess(g):g< f\}\geq \frac{ess(f)}{3}.
\end{displaymath}
\end{enumerate}
\end{lemma}

\begin{proof} The first statement follows immediately from the definition of $\leq  $.
To see that the statement $2.$ also holds, let $f$ be an $n$-ary Boolean function $f$ with $ess(f)\geq3$.
Without loss of generality, we may suppose that $n=ess(f)$. For each $i,j=1,2,3$, $i\not=j$, let $f_{ij}$ be the function obtained from $f$ by identifying the $i$th and
$j$th variables of $f$. To avoid notational difficulties, we will not relabel the variables. 
We claim that $ess(f_{ij})\geq \frac{n}{3}$ for some pair $ij$.
{\it Case 1.} $n=3$. If  $ess(f_{ij})\geq 1$ for no pair $ij$ then $f_{12}, f_{13}, f_{23}$ are constant. In fact, they take the same value $a:= f_{123}$. But, since the first variable of $f$ is essential, we have $f(0, a_2,a_3)\not =f(1, a_2,a_3)$ for some $a_2,a_3\in \B$. Since $\B$ has two elements,  each of the triples $(0, a_2,a_3)$, $(1, a_2,a_3)$ has  two components which are equal, thus $f(0, a_2,a_3)=f(1, a_2,a_3)=a$, a contradiction. 
{\it Case 2.} $n\geq 4$. Let $A_{ij}:=\{k:   4\leq k\leq n\; \text{and}\; x_k \; \text{is an essential variable of }\; f_{ij}\}$. If  for some $i,j=1,2,3$, $i\not=j$, $A_{ij}$ has at least $\frac{n}{3}$ elements, the claimed inequality is proved. If not, we claim that all $A_{ij}$ have $\frac{n}{3}-1$ elements. For that it suffices to observe that  each $k$, $4\leq k\leq n$,  belongs to some $A_{ij}$ for some
$i,j=1,2,3$, $i\not=j$.   This fact is easy to obtain. Since $x_k$ is  an essential variable of $f$ there are 
$(a_1, a_2, a_3,\ldots ,a_{k-1} , a_{k+1},\dots , a_n)$ such that 
  \begin{displaymath}\label{eqessential}
 f(a_1,a_2,a_3,...,a_{k-1},0 ,a_{k+1},\dots, a_n)\not = f(a_1,a_1,a_3,\dots,a_{k-1},1 ,a_{k+1}, ..., a_n)
\end{displaymath}
Since $\mathbb{B}$ has two elements, there are $i,j$ with $1\leq i<j\leq 3$ such that 
$a_i=a_j$ for some pair $\{i,j\}$. Therefore, there are 
 $(b_1, b_2, b_3,\ldots ,b_{k-2} , b_{k},\ldots , b_{n-1} )$ such that 
  \begin{displaymath}
 f_{ij}(b_1,b_2,\ldots ,b_{k-2},0 ,b_{k},\dots, b_{n-1}) \not = f_{ij}(b_1,b_2,\ldots ,b_{k-2},1 ,b_{k},\dots, b_{n-1}) 
\end{displaymath}
Hence, $x_k$ (which corresponds to the $k-1$th variable of $f_{ij}$) is an essential variable of $f_{ij}$,
proving that our observation holds. 
 Thus, for every pair $i,j=1,2,3$, $i\not=j$ there are $\frac {n} {3}-1$ variables $k$ which
 are essential   for $f_{ij}$. As in the proof of Case 1, there is some pair $ij$, where $1\leq i,j \leq 3$,
 for which $1\leq l\leq 3$ is also essential for $f_{ij}$. Hence, 
 $f_{ij}$ has $\frac{n}{3}$ essential variables as claimed. 
\end{proof}

 \begin{corollary}\label{levels} In $(\Omega /\equiv , \sqsubseteq )$ every principal initial 
segment is finite and each level is finite.   
\end{corollary}
\begin{proof}
According to the above lemma, for every $n\geq 1$, and for each Boolean function
$f$ in the $n$-th level, we have $n<ess(f)\leq 3^n$. The result follows. 
\end{proof}

\subsection{Definability of Boolean function classes by means of functional equations}

A \emph{functional equation} (for Boolean functions) is a formal expression
\begin{equation}\label{functionalEq}
\begin{array}{l}
h_1( {\bf f} (g_1({\bf x}_1,\ldots ,{\bf x}_p)),\ldots ,{\bf f} (g_m({\bf x}_1,\ldots ,{ \bf x}_p)))=
\\=h_2( {\bf f} ({g'}_1({\bf x}_1,\ldots ,{\bf x}_p)),\ldots ,{\bf f} ({g'}_t({\bf x}_1,\ldots ,{\bf x}_p))) 
\end{array}
\end{equation}
where $m,t,p\geq 1$, $h_1: \mathbb{B}^m \rightarrow \mathbb{B}$,
$h_2: \mathbb{B}^t \rightarrow \mathbb{B}$, each $g_i$ and ${g'}_j$ is a map
$\mathbb{B}^p\rightarrow \mathbb{B}$, the ${\bf x}_1,\ldots ,{\bf x}_p$ are $p$ distinct
\emph{vector variable symbols}, and $\bf f$ is a distinct \emph{function symbol}.
Such equations were systematically studied in \cite{EFHH}. See e.g. \cite{Po,FP,Pi2} for variants, and
\cite{CF2} for extensions and more stringent notions of functional equations.

An $n$-ary Boolean function $f: \mathbb{B}^n \rightarrow \mathbb{B}$,
\emph{satisfies} the equation (\ref{functionalEq}) if, for all ${\bf v}_1,\ldots ,{\bf v}_p\in \mathbb{B}^n$, we have
 \begin{displaymath}
\begin{array}{l}
h_1( {f} (g_1({\bf v}_1,\ldots ,{\bf v}_p)),\ldots ,{f} (g_m({\bf v}_1,\ldots ,{\bf v}_p)))=
\\=h_2( {f} ({g'}_1({\bf v}_1,\ldots ,{\bf v}_p)),\ldots ,{f} ({g'}_t({\bf v}_1,\ldots ,{\bf v}_p)))
\end{array}
\end{displaymath}
where $g_1({\bf v}_1,\ldots ,{\bf v}_p)$ is interpreted component-wise, that is,
 \begin{displaymath}
\begin{array}{l}
g_1({\bf v}_1,\ldots ,{\bf v}_p)=(g_1({\bf v}_1(1),\ldots ,{\bf v}_p(1)),\ldots ,g_1({\bf v}_1(n),\ldots ,{\bf v}_p(n)))
\end{array}
\end{displaymath}
 A class $K$ of Boolean functions is said to be \emph{defined}
 by a set $\mathcal{E}$ of functional equations, if $K$ is the class of all
those Boolean functions which satisfy every member of $\mathcal{E}$.
It is not difficult to see that if a class $K$ is defined by a set $\mathcal{E}$ of functional equations,
 then it is also defined by a set $\mathcal{E}'$ whose members are functional equations in which the
indices $m$ and $t$ are the same.

By an \emph{equational class} we simply mean a class of Boolean functions definable by a set
of functional equations.
The following characterization of equational classes
was first obtained by Ekin, Foldes, Hammer and Hellerstein \cite{EFHH}. For variants and extensions,
 see e.g.\cite{FP,Po,CF2}.

\begin{theorem} The equational classes of
Boolean functions are exactly those classes
 that are closed under simple variable substitutions.
\end{theorem}
In other words, a class $K$ is equational if and only if $K/\equiv $ is an initial segment of
$\Omega /\equiv $.

\subsection{Definability of Boolean function classes by means of relational constraints}

An \emph{$m$-ary Boolean relation} is a subset $R$ of $ \mathbb{B}^m$.
Let $f$ be an $n$-ary Boolean function. We denote by $fR$ the $m$-ary relation given by
 \begin{displaymath}
\begin{array}{l}
fR=\{ f( {\bf v}_1,\ldots ,{\bf v}_n) : {\bf v}_1,\ldots ,{\bf v}_n\in R\}
\end{array}
\end{displaymath}
 where the $m$-vector $f( {\bf v}_1,\ldots ,{\bf v}_n) $ is defined component-wise as
in the previous subsection.

An $m$-ary \emph{Boolean constraint}, or simply an $m$-ary \emph{constraint},
is a pair $(R,S)$ where $R$ and $S$ are $m$-ary relations called the \emph{antecedent} and
\emph{consequent}, respectively,
 of the relational constraint.
A Boolean function is said to \emph{satisfy}
 an $m$-ary constraint $(R,S)$ if $fR\subseteq S$.
 Within this framework, a class $K$ of Boolean functions is said to be \emph{defined}
 by a set $\mathcal{T}$ of relational constraints, if $K$ is the class of all
those Boolean functions which satisfy every member of $\mathcal{T}$.
For further background, see \cite{Pi2}.
See also \cite{C2,CF1, CF2, CF3, L}, for further variants and extensions.

The connection between definability by functional equations and by relational constraints was made explicit by Pippenger who
established in \cite{Pi2} a complete correspondence between functional equations and relational constraints.
\begin{theorem}
The equational classes of Boolean functions are exactly those classes definable by relational constraints.
\end{theorem}
 This result was further extended and strengthened in \cite{CF3}.

\begin{proposition}\label{prop} For each relational constraint $(R,S)$ there is a functional equation satisfied by
 exactly the same Boolean functions satisfying $(R,S)$. Conversely, for each functional equational
 \begin{displaymath}
\begin{array}{l}
h_1( {\bf f} (g_1({\bf x}_1,\ldots ,{\bf x}_p)),\ldots ,{\bf f} (g_m({\bf x}_1,\ldots ,{ \bf x}_p)))=
\\=h_2( {\bf f} ({g'}_1({\bf x}_1,\ldots ,{\bf x}_p)),\ldots ,{\bf f} ({g'}_t({\bf x}_1,\ldots ,{\bf x}_p))) \qquad \qquad (\ref{functionalEq})
\end{array}
\end{displaymath}
there is a relational constraint satisfied by exactly the same Boolean functions satisfying (\ref{functionalEq}).
\end{proposition}

\begin{proof} We follow the same steps as in the proof of Theorem 1 in \cite{CF3}.
For each functional equation (\ref{functionalEq}),
 let $(R,S)$ be the relational constraint defined by
\begin{displaymath}
\begin{array}{l}
R:=\{ (g_1({\bf a}),\ldots ,g_m({\bf a}),{g'}_1({\bf a}),\ldots ,{g'}_t({\bf a})) :{\bf a}\in \mathbb{B}^p\},
\\S:=\{(b_1,\ldots ,b_m,{b'}_1,\ldots ,{b'}_t)\in \mathbb{B}^{m+t}: h_1(b_1,\ldots ,b_m)=h_2({b'}_1,\ldots ,{b'}_t)\}.
\end{array}
\end{displaymath}
Let $f$ be an $n$-ary Boolean function. From the definition of $S$, it follows that 
$f$ satisfies $(R,S)$ if and only if for every ${\bf a}_1,\ldots ,{\bf a}_n\in R$,
 \begin{displaymath}
\begin{array}{l}
h_1( {f} ({\bf a}_1(1),\ldots ,{\bf a}_n(1)),\ldots ,{f} ({\bf a}_1(m),\ldots ,{\bf a}_n(m)))=
\\=h_2( {f} ({\bf a}_1(m+1),\ldots ,{\bf a}_n(m+1)),\ldots ,{f} ({\bf a}_1(m+t),\ldots ,{\bf a}_n(m+t))) 
\end{array}
\end{displaymath}
Since $R$ is the range of $g=(g_1,\ldots ,g_m,{g'}_1,\ldots ,{g'}_t)$, we have that
$f$ satisfies $(R,S)$ if and only if for every ${\bf v}_1,\ldots ,{\bf v}_p\in \mathbb{B}^n$
 \begin{displaymath}
\begin{array}{l}
h_1( {f} (g_1({\bf v}_1,\ldots ,{\bf v}_p)),\ldots ,{f} (g_m({\bf v}_1,\ldots ,{\bf v}_p)))=
\\=h_2( {f} ({g'}_1({\bf v}_1,\ldots ,{\bf v}_p)),\ldots ,{f} ({g'}_t({\bf v}_1,\ldots ,{\bf v}_p)))
\end{array}
\end{displaymath}
In other words, $f$ satisfies $(R,S)$ if and only if $f$ satisfies (\ref{functionalEq}).

Conversely, let 
$(R,S)$ be a relational constraint. We may suppose $R$ non-empty, indeed,
 constraints with empty antecedent are satisfied by every 
Boolean function, and thus they can be discarded as irrelevant.  
With the help of the following two facts, we will construct a functional equation satisfied by the exactly the
same functions as those satisfying $(R,S)$.
\begin{fact} For each non-empty Boolean relation $R\subseteq \mathbb{B}^m$, there is a $p\geq 1$ and a map
$g:= (g_1,\ldots ,g_m)$, where each $g_i$ is a $p$-ary Boolean function
$g_i:\mathbb{B}^p\rightarrow \mathbb{B}$, such that the range of $g$ is $R$.
\end{fact}

\begin{fact} For each Boolean relation $S\subseteq \mathbb{B}^m$, there exist maps
$h_1,h_2: \mathbb{B}^m\rightarrow \mathbb{B}$, such that
\begin{displaymath}
S=\{{\bf b}\in B^{m}: h_1({\bf b})=h_2({\bf b})\}.
\end{displaymath}
\end{fact}

Let $(R,S)$ be a relational constraint. Consider the functional equation
 \begin{equation}\label{constructed}
\begin{array}{l}
h_1( {\bf f} (g_1({\bf x}_1,\ldots ,{\bf x}_p)),\ldots ,{\bf f} (g_m({\bf x}_1,\ldots ,{ \bf x}_p)))=
\\=h_2( {\bf f} ({g}_1({\bf x}_1,\ldots ,{\bf x}_p)),\ldots ,{\bf f} ({g}_m({\bf x}_1,\ldots ,{\bf x}_p))) 
\end{array}
\end{equation}
where the $g_i$'s and $h_j$'s are the maps given in Fact 1 and Fact 2.
Let $f$ be an $n$-ary Boolean function. By construction, we have that $f$ satisfies 
(\ref{constructed}) if and only if for every ${\bf v}_1,\ldots ,{\bf v}_p\in \mathbb{B}^n$,
 $(f (g_1({\bf v}_1,\ldots ,{\bf v}_p)),\ldots ,{ f} (g_m({\bf v}_1,\ldots ,{ \bf v}_p)))\in S$.
>From the fact that $R$ is the range of $(g_1,\ldots ,g_m)$, it follows that 
 $f$ satisfies (\ref{constructed}) if and only if  $f$ satisfies $(R,S)$.
\end{proof}

In the sequel, we will make use of the following result  of Pippenger (\cite{Pi2}, Theorem 2.1). For the reader convenience, we provide a proof. 
\begin{lemma}\label{pippenger} For each Boolean function $f$,
there is a relational constraint $(R,S)$ such that $\Omega(R,S)= Forbid(\{f\})$.
 \end{lemma}
 \begin{proof}
Let $f$ be Boolean function, say of arity $n$. Let
${\bf v}_1,\ldots ,{\bf v}_n$ be $2^n$-vectors such that
$\mathbb{B}^n= \{({\bf v}_1(i),\ldots ,{\bf v}_n(i)): 1\leq i\leq 2^n\}$. Consider the $2^n$-ary relations
$R_f$ and $S_f$ given by
\begin{displaymath}
\begin{array}{l}
R: =\{ {\bf v}_1,\ldots ,{\bf v}_n\},\textrm{ and } S_f:=  \bigcup \{g R_f :  g\in Forbid(\{f\})\}
\end{array}
\end{displaymath}
respectively. Clearly, if $g\in Forbid(\{f\})$, 
then $g$ satisfies $(R_f,S_f)$.
If $g'$, say $m$-ary, is a member of $\uparrow f$, then there are $n$-ary projections
$p_1, \ldots ,p_m\in I_c$ such that  
\begin{equation}\label{equationminor}
f=g'(p_1,\ldots ,p_m)
\end{equation}
We claim that $g'(p_1({\bf v}_1,\ldots , {\bf v}_n),\dots, p_m({\bf v}_1,\ldots ,{\bf v}_n))$ does not belong to $S_f$.
 Otherwise, there would be $g \in Forbid(\{f\})$, and projections $p'_1,\ldots , p'_t$ such that
 \begin{displaymath}
\begin{array}{l}
 g'(p_1({\bf v}_1,\ldots , {\bf v}_n),\ldots , p_m({\bf v}_1,\ldots ,{\bf v}_n))=\\
g(p'_1({\bf v}_1,\ldots , {\bf v}_n),\ldots ,  p'_t({\bf v}_1,\ldots ,{\bf v}_n)).
 \end{array}
\end{displaymath}
By definition, this amounts to  
 \begin{displaymath}
\begin{array}{l}
g'(p_1({\bf v}_1,\ldots , {\bf v}_n)(i),\ldots , p_m({\bf v}_1,\ldots , {\bf v}_n)(i))=\\
g(p'_1({\bf v}_1,\ldots , {\bf v}_n)(i),\ldots , p'_t({\bf v}_1,\ldots , {\bf v}_n)(i))
  \end{array}
\end{displaymath}
for all $i$, $1\leq i\leq 2^n$. Which, in turn, amounts to 
 \begin{displaymath}
\begin{array}{l}
g'(p_1({\bf v}_1(i), \ldots , {\bf v}_n(i)), \ldots , p_m({\bf v}_1(i),\ldots , {\bf v}_n(i)))=\\
g(p'_1({\bf v}_1(i),\ldots , {\bf v}_n(i)),\ldots , p'_t({\bf v}_1(i),\ldots ,{\bf v}_n(i)))
 \end{array}
\end{displaymath} 
Since for every $(x_1,\ldots , x_n)\in \B^n$  there is some $i$ such that 
 \begin{displaymath}
\begin{array}{l}
({\bf v}_1(i),\ldots , {\bf v}_n(i))= (x_1,\ldots , x_n)
 \end{array}
\end{displaymath} 
 we get 
 \begin{displaymath}
\begin{array}{l}
g'(p_1(x_1,\ldots , x_n),\ldots , p_m(x_1,\ldots ,x_n))=\\
g(p'_1(x_1,\ldots , x_n),\ldots , p'_t(x_1,\ldots ,x_n))
\end{array} 
\end{displaymath} 
 for all $(x_1,\ldots ,x_n)\in \B^n$, that is 
\begin{displaymath}\label{identg}
g'(p_1,\ldots , p_m)=g(p'_1,\ldots , p'_t)
\end{displaymath} 
With equation (\ref{equationminor}) we get 
$f= g(p'_1,\ldots , p'_t)$ that is $f$ is 
obtained from $g$ by simple variable substitutions, contradicting our assumption 
$g \in Forbid(\{f\})$. 
\end{proof}

\section{Proofs}
\subsection{Proof of Theorem \ref{thm2}}
We show that $(i)\Rightarrow (ii)\Rightarrow (iii)\Rightarrow (i)$.

$(i)\Rightarrow (ii)$ 
To see that each class $I\in \mathcal{F}$ can be defined by a single functional equation, note that
 \begin{displaymath}
\begin{array}{l}
h_1( {\bf f} (g_1({\bf x}_1,\ldots ,{\bf x}_p)),\ldots ,{\bf f} (g_m({\bf x}_1,\ldots ,{ \bf x}_p)))=
\\=h_2( {\bf f} ({g'}_1({\bf x}_1,\ldots ,{\bf x}_p)),\ldots ,{\bf f} ({g'}_t({\bf x}_1,\ldots ,{\bf x}_p))) \qquad (1)
\end{array}
\end{displaymath}
is satisfied by exactly the same functions satisfying
 \begin{displaymath}
\begin{array}{l}
h_1( {\bf f} (g_1({\bf x}_1,\ldots ,{\bf x}_p)),\ldots ,{\bf f} (g_m({\bf x}_1,\ldots ,{ \bf x}_p)))+ \\
h_2( {\bf f} ({g'}_1({\bf x}_1,\ldots ,{\bf x}_p)),\ldots ,{\bf f} ({g'}_t({\bf x}_1,\ldots ,{\bf x}_p)))=0
\end{array}
\end{displaymath}
where $+$ denotes the sum modulo 2. Thus, if $I$ is defined by the equations $H_1=0,\ldots ,H_n=0$,
 then it is also defined by $\bigvee _{1\leq i\leq n}H_i=0$.

 $(ii)\Rightarrow (iii)$ Let $L$ be a functional equation. According to Proposition \ref{prop}, there is a relational constraint $(R,S)$ such that the operations satisfying   $\Omega(R,S)$ are those satisfying 
 $L$.
 \begin{lemma}The set $\Omega(R,S)$ of operations  which satisfy a $n$-ary constraint $(R,S)$ is of the form $Forbid(A)$ for some finite antichain  $A$ of $\Omega$.
\end{lemma}
\begin{proof}
\begin{claim}\label{claim1}
 If an $m$-ary Boolean function $g$ does not satisfy $(R, S)$, 
 then there is some $m'$-ary $g'$, where $m'\leq 2^n$, such that $g'\leq g$
and such that $g'$ does not satisfy $(R, S)$.
\end{claim}
\begin{proof}[Proof of Claim \ref{claim1}] If $m\leq 2^n$ set $g':= g$. If not, let  $v_1, \ldots , v_m\in R$ such that
 $g(v_1, \ldots , v_m)\not \in S$.
 Say that two indices $i,j$ with $1\leq i,j\leq m$ are equivalent if  
 $v_i=v_j$.  Let $C_1,\ldots , C_{m'}$ be an enumeration of the equivalence  classes. 
For each $i$, $1\leq i\leq m$,  let $c(i)$ be the indice for which $i\in C_c(i)$. 
Let $g'$ be the $m'$-ary operation defined by $g':= g(p_1, \ldots , p_m)$, where 
$p_j (x_1,\ldots , x_{m'})=x_c(j)$. Clearly, $m'\leq 2^n$ and, by definition, $g'\leq g$. 
For each $1\leq j\leq m'$, let $w_j:= v_i$, whenever $c(i)=j$. We have 
$g(v_1, \ldots , v_m)= g(w_{c(1)}, \ldots , w_{c(m)})$ and since  
$g'(x_1,\ldots ,x_{m'})= g(x_c(1),\ldots , x_c(m))$ it follows that
$g'(w_1, \ldots , w_m)= g(v_1, \ldots , v_m)$ and hence, $g'$ does not satisfy $(R,S)$.
\end{proof}
From Claim \ref{claim1}, the minimal members of $\Omega\setminus \Omega(R,S)$ 
 have arity at most $2^n$ and hence,
 there are only finitely many of such minimal members (w.r.t. the equivalence associated with the quasi-order).
\end{proof}

$(iii)\Rightarrow (i)$ Let $I:= Forbid (A)$ where $A$ is a finite antichain. Since  $I$ is  a finite intersection of set of the form $Forbid(\{f\})$, in order to get that $I\in \mathcal F$, it suffices to show that $Forbid(\{f\})\in \mathcal F$. 
This is a consequence of Proposition \ref{prop} and Lemma \ref{pippenger}.

\subsection{Proof of Lemma \ref{mainlemma}}

Statement $(1)$. 
If $K_1$ and $K_2$ are classes in $\mathcal{F}$, say
defined by the expressions
 \begin{displaymath}
\begin{array}{l}
H_1=0 \textrm{ and } H_2=0
\end{array}
\end{displaymath}
respectively, then $K_1\cup K_2$ and $K_1\cap K_2$ are defined by
 \begin{displaymath}
\begin{array}{l}
H_1\wedge H_2=0 \textrm{ and } H_1\vee H_2=0
\end{array}
\end{displaymath}
respectively. This proves that statement $(1)$ of Lemma \ref {mainlemma} holds. 
 The fact that $\mathcal F$ is closed under finite intersections follows also from 
 the equivalence $(i)\Rightarrow (iii)$ of  Theorem \ref{thm2}.
 Note that from this equivalence and the fact that $\mathcal F$ is closed under finite unions, it  follows that $\Omega/\equiv $ is up-closed.
 
 Statement $(2)$ Implication $(iii)\Rightarrow (i)$ of Theorem \ref{thm2}.  
 
Statement $(3)$ Let $f\in \Omega$.  Let  $\overline f$ be its image in 
$P:= (\Omega/\equiv, \sqsubseteq)$ (i.e., the equivalence class containing $f$), and $m:=h(\overline f, P)$.  The initial segment $\downarrow f$ is of the form $Forbid (A)$ for some antichain $A$. This antichain $A$ is made of representative of the minimal elements of $B:= P\setminus \downarrow \overline f$. If $y$ is minimal in  $B$ then for every $x$ such that $x<y$, we have $x\leq \overline f$. It follows that $h(y, P)\leq h(\overline f, P)+1=m+1$, that is the minimal elements of $B$ belong to  the union of levels $P_n$ for $n\leq m+1$. From Corollary \ref{levels}, all levels of $P$ are finite. Hence $A$ is finite.

Statement $(4)$ Let $E^k$ be the set of operations with at most $k$ essential variables. Its image $\overline E^k$ in $P:= (\Omega/\equiv, \sqsubseteq)$ is in fact included into the union of all levels $P_n$ for $n\leq k$. Since by  Corollary \ref{levels}, all levels are finite,  $E^k$ is a finite union of initial segments of the form $\downarrow f$. According to  Statement $(1)$ and Statement $(3)$, $E^k \in \mathcal F$. 
\endproof

\section{Proof of Theorem \ref{main}}
Let $P:= (\Omega/\equiv, \sqsubseteq)$.

Part 1. $P$ embeds into $([\omega]^{<\omega}, \subseteq)$.

We	apply Lemma \ref{finite}. The poset $P$ is trivially countable, and by Corollary \ref{levels}, 
 for every $x\in P$, the initial segment $\downarrow x$ is finite. Thus, by Lemma \ref{finite}, 
$P$ embeds into $([\omega]^{<\omega}, \subseteq)$.

Part 2. $([\omega]^{<\omega}, \subseteq)$ embeds into $P$. The following is a particular case of Proposition 3.4 in \cite{Pi2}.
\begin{lemma}\label{lemma1}
The family $(f_n)_{n\geq 4}$ of Boolean functions, given by
\begin{displaymath}
f_n(x_1,\ldots ,x_n)=\left\{
             \begin{array}{ll}
                     1         &       \mbox{if $\#\{i:x_i=1\}\in \{1,n-1\}$}   \\
                     0         &       \mbox{otherwise}.
             \end{array}\right.
 \end{displaymath}
constitutes an  infinite antichain of Boolean functions.
\end{lemma}

Note that $f_n(a,\ldots ,a)=0$ for $a\in \{0,1\}$. The following lemma was presented in \cite{C3}.

\begin{lemma}\label{lemma3}
Let $(f_n)_{n\geq 4}$ be the family of Boolean functions
given above, and consider the family $(u_n)_{n\geq 4}$ defined by
 \begin{displaymath}
\begin{array}{llll}
    u_n(x_0,x_1,\ldots ,x_n) =x_0\wedge f_n(x_1,\ldots ,x_n)     \\
 \end{array}
 \end{displaymath}
 The family $(u_n)_{n\geq 4}$ constitutes an infinite antichain of Boolean functions.
\end{lemma}

\begin{proof} We follow the same steps as in \cite{C3}. We show that if $m\not=n$, then
$u_m\not\leq u_n$.
 By definition, $u_m$ and $u_n$ cannot have dummy variables. Therefore,
 $u_m\not\leq u_n$, whenever $m>n$.

 So assume that $m<n$, and for a contradiction, suppose that $u_m\leq u_n$, i.e.
there are $m+1$-ary projections $p_0,\ldots ,p_n\in I_c$ such that
$u_m=u_n(p_0,\ldots ,p_n)$. Note that for every $m\geq 4$,
$u_m(1,x_1\ldots ,x_m)=f_m(x_1\ldots ,x_m)$ and $u_m(0,x_1\ldots ,x_m)$ is the constant 0.

Now, suppose that $p_0(x_0,\ldots ,x_m)=x_0$. If for all $1\leq k\leq n$, $p_k(x_0,\ldots ,x_m)\not=x_0$,
then by taking $x_0=1$ we would conclude that $f_m\leq f_n$,
contradicting Lemma \ref{lemma1}. If there is $1\leq k\leq n$ such that $p_k(x_0,\ldots ,x_m)=x_0$, then
  by taking $x_i=1$ if and only if $i=0,1$, we would have
 \begin{displaymath}
 \begin{array}{ll}
u_m(x_0,\ldots  ,x_m)=1\not=0=
u_n(x_0,p_1(x_0,\ldots  ,x_m),\ldots ,p_n(x_0,\ldots  ,x_m))
 \end{array}
 \end{displaymath}
which is also a contradiction.

Hence, $p_0(x_0,\ldots ,x_m)\not=x_0$, say $p_0(x_0,\ldots ,x_m)=x_i$ for $1\leq i\leq m$. But then
 by taking $x_i=1$ if and only if $i=0,1$, we would have
\begin{displaymath}
 \begin{array}{ll}
u_m(x_0,\ldots ,x_i,\ldots  ,x_m)=1\not=\\
0=u_n(x_i,p_1(x_0,\ldots ,\ldots  ,x_m),\ldots ,p_n(x_0,\ldots ,x_m))=0
 \end{array}
 \end{displaymath}
 which contradicts our assumption $u_m\leq u_n$.
\end{proof}

Let $I$ be a non-empty finite set of integers greater or equal than 4, and let $g_I$ be the
$ \sum _{i\in I}i$-ary function given by
  \begin{displaymath}
 \begin{array}{ll}
g_I= \sum_{i\in I} \bigwedge_{j\in I\setminus \{i\}}\bigwedge_{1\leq k\leq j} x_k^j \wedge f_i(x_1^i,\ldots ,x_i^i)
  \end{array}
\end{displaymath}
Observe that
\begin{itemize}
\item By identifying all $x_k^j$, for $j\in I\setminus \{i\}$ and $1\leq k\leq j$, we obtain
$x_0 \wedge f_i(x_1^i,\ldots ,x_i^i)$, and
\item $g_I=1$ if and only if there exactly one $i\in I$ such that
\item[$i)$] for all $j\in I\setminus \{i\}$ and $1\leq k\leq j$, $x_k^j =1$, and
\item[$ii)$] $\#\{1\leq k\leq i:x_k^i=1\}\in \{1,i-1\}$.
\end{itemize}

\begin{proposition} \label{mainProp}
Let $I$ be a non-empty finite set of integers greater or equal than 4,
and let $g_I$ be the $ \sum _{i\in I}i$-ary function given above. Then for every $n\geq 4$,
$n\in I$ if and only if $u_n\leq g_I$.
\end{proposition}

\begin{proof} By the first observation above it follows that if $n\in I$ then $u_n\leq g_I$.
To prove the converse, suppose that $n\not\in I$ and for a contradiction suppose that
$u_n\leq  g_I$, i.e., there are projections $p_k^i(x_0,x_1,\ldots ,x_n)$,
 $i\in I$ and $1\leq k\leq i$, such that
 \begin{equation}\label{construction}
\begin{array}{l}
u_n=\sum_{i\in I}\bigwedge_{j\in I\setminus \{i\}}\bigwedge_{1\leq k\leq j} p_k^j \wedge f_i(p_1^i,\ldots ,p_i^i) 
\end{array}
\end{equation}

Consider the vector $(a_0,a_1, \ldots ,a_n)$ given by $a_l=1$ iff $l=0,1$.

Clearly, $u_n(a_0,a_1,\ldots  ,a_n)=1$ and, in order to have $(\ref{construction})=1$,
 there must exist exactly one $i\in I$ such that
\begin{itemize}
\item[$i)$] for all $j\in I\setminus \{i\}$ and $1\leq k\leq j$, $p_k^j \in \{x_0,x_1\}$, and
\item[$ii)$] $\#\{1\leq k\leq i:p_k^i \in \{x_0,x_1\}\}\in \{1,i-1\}$.
\end{itemize}
Since $x_2,\ldots ,x_n$ are essential in $u_n$, we also have that for each $2\leq l\leq n$, there is
$1\leq k\leq i$ such that $p_k^i=x_l$.

Now, if for all $1\leq k\leq i$, $p_k^i \not=x_1$, then there are
$j\in I\setminus \{i\}$ and $1\leq k\leq j$, such that $p_k^j =x_1$, because $x_1$ is essential in $u_n$.
Consider $(b_0,b_1, \ldots ,b_n)$
given by $b_l=1$ iff $l=0,2$. We have $u_n(b_0,b_1,\ldots  ,b_n)=1$,
 but $(\ref{construction})=0$, which constitutes a contradiction.

Thus, there is $1\leq k\leq i$, $p_k^i =x_1$. If there is $1\leq t\leq i$ such that $p_t^i =x_0$, then for 
$(b_0,b_1, \ldots ,b_n)$ given by $b_l=1$ iff $l=0,1$, 
we have $u_n(b_0,b_1,\ldots  ,b_n)=1$,
 but $(\ref{construction})=0$, because for each $2\leq l\leq n$, there is
$1\leq r\leq i$ such that $p_r^i=x_l$ and $n\geq 4$. Hence, for every $1\leq t\leq i$, $p_t^i \not=x_0$, 
and since for each $2\leq l\leq n$, there is $1\leq r\leq i$ such that $p_r^i=x_l$, we must have
$i\not<n$. Also, $n\not\in I$ and thus $i>n$. But in this case, there must exist $1\leq s\leq n$ such that,
for some $1\leq r_1<r_2\leq i$, $p_{r_1}^i=p_{r_2}^i=x_s$. Now, if $s=1$, then for 
$(b_0,b_1, \ldots ,b_n)$ given by $b_l=1$ iff $l=0,1$, 
we have $u_n(b_0,b_1,\ldots  ,b_n)=1$ and $(\ref{construction})=0$, once again by the fact that
 for each $2\leq l\leq n$, there is $1\leq r\leq i$ such that $p_r^i=x_l$.
If $1< s\leq n$, then for $(b_0,b_1, \ldots ,b_n)$ given by $b_l=1$ iff $l\not=s$, 
we have $u_n(b_0,b_1,\ldots  ,b_n)=1$
 and $(\ref{construction})=0$.

Since in all possible cases we derive the same contradiction,
the proof of the proposition is complete.
\end{proof}

By making use of Proposition \ref{mainProp}, it is not difficult to verify that the mapping
$I\mapsto g_{I'}$, where
$I'=\{i+4:i\in I\}$, is an embedding from  $( [\omega ]^{<\omega },\subseteq )$ into $(\Omega /\equiv , \leq  )$.

\section{Linear functions with a bounded number of essential variables}
\begin{theorem}\label{characterization} The class $L^k$ of linear functions with at most
$k\geq 1$ essential variables is defined by
\begin{equation}\label{EqCharacterizing}
\begin{array}{l}
( \bar {\bf f}({\bf 0}) \wedge (\bigwedge_{1\leq i\leq k+1}{\bf f}({\bf x}_i)\longrightarrow \bigvee _{1\leq j<l\leq k+1}{\bf f}({\bf x}_j{\bf x}_l)) )\quad \vee
\\  ({\bf f}({\bf 0}) \wedge (\bigwedge_{1\leq i\leq k+1}({\bf f}({\bf x}_i)+1)\longrightarrow 
\bigvee _{1\leq j<l\leq k+1}({\bf f}({\bf x}_j{\bf x}_l)+1)) )=1 
\end{array}
\end{equation}
\end{theorem}

\begin{proof}
Note that $L^k$ is the class of linear functions which are the sum of at most $k\geq 0$ variables.
First we show that if $f\in L\setminus L^k$, then $f$ does not satisfy (\ref{EqCharacterizing}).
So suppose that $f$ is the sum of $k\geq 1$ variables. Without loss of generality, assume that
$f=x_1+\ldots +x_{k+1}+c_{k+2}x_{k+2}+\ldots +c_nx_n+c$, where $c_{k+2},\ldots ,c_{n},c\in \{0,1\}$.
 For $1\leq i\leq k+1$, let ${\bf a}_i$ be the unit $n$-vector
with all but the $i$-th component equal to $0$. If $c=0$, then $ { f}({\bf 0})=0$ and
 \begin{displaymath}
\begin{array}{l}
(\bigwedge_{1\leq i\leq k+1}{ f}({\bf a}_i)\longrightarrow \bigvee _{1\leq j<l\leq k+1}{ f}({\bf a}_j{\bf a}_l))=0
\end{array}
\end{displaymath}
Thus $f$ does not satisfy $(3)$. If $c=1$, then $ {f}({\bf 0})=1$ and
 \begin{displaymath}
\begin{array}{l}
(\bigwedge_{1\leq i\leq k+1}({ f}({\bf a}_i)+1)\longrightarrow \bigvee _{1\leq j<l\leq k+1}({ f}({\bf a}_j{\bf a}_l)+1)=0
\end{array}
\end{displaymath}
Thus $f$ does not satisfy (\ref{EqCharacterizing}).

Now we show that every linear function $f$ in $L^k$ satisfies (\ref{EqCharacterizing}).
We make use of the following
 \begin{claim}\label{claim2}
Let $1\leq n\leq k$ and let ${\bf x}_1,\ldots ,{\bf x}_{k+1}$ be $k+1 $ $n$-vectors of odd weight.
Then there are $1\leq i,j\leq k+1$, $i\not=j$, 
such that ${\bf x}_j{\bf x}_i$ has odd weight. 
\end{claim}
\begin{proof}[Proof of Claim \ref{claim2}.]
Let ${\bf x}_1,\ldots ,{\bf x}_{k+1}$ be $k+1 $ $n$-vectors of odd weight.
Since there are at most $n$
linearly independent $n$-vectors, ${\bf x}_1,\ldots ,{\bf x}_{k+1}$ must be linearly dependent,
 i.e., there is $I\subseteq \{1,\ldots ,k+1\}$ and
$ j\in \{1,\ldots ,k+1\}\setminus I$ such that ${\bf x}_j=\sum _{i\in I}{\bf x}_i$. 
We have 
${\bf x}_j={\bf x}_j {\bf x}_j={\bf x}_j\sum _{i\in I}{\bf x}_i= \sum _{i\in I}{\bf x}_j{\bf x}_i $.
 Since the weight of ${\bf x}_j$ is odd, and the weight function modulo $2$ (i.e. the parity function)
 distributes over the component-wise sum of vectors, it follows that there 
is an odd number of products ${\bf x}_j{\bf x}_i$,
 $i\in I$, with odd weight. In particular, there are $1\leq i,j\leq k+1$, $i\not=j$, 
such that ${\bf x}_j{\bf x}_i$ has odd weight. 
\end{proof}

Let $f $ be a linear function in $L^k$, say $f=x_1+\ldots +x_{n}+c$,
where $c\in \{0,1\}$ and $1\leq n\leq k$. If $c=0$, then $ { f}({\bf 0})=0$ and
$ { f}( {\bf a})=1$ if and only if ${\bf a}$ has odd weight.
Now, if ${\bf a}_1,\ldots ,{\bf a}_{k+1}$ are $k+1 $ $n$-vectores such that
$\bigwedge_{1\leq i\leq k+1}{ f}({\bf a}_i)=1$, then each ${\bf a}_i$, $1\leq i\leq k+1$,
has odd weight and by Claim \ref{claim2} it follows that there are
$1\leq i<j\leq k+1$ such that ${\bf a}_i{\bf a}_j$ has odd weight, and hence,
$\bigvee _{1\leq j<l\leq k+1}{ f}({\bf a}_j{\bf a}_l)=1$.

If $c=1$, then $ { f}({\bf 0})=1$ and
$ { f}( {\bf a})+1=1$ if and only if ${\bf a}$ has odd weight. Again, by making use of Claim \ref{claim2},
it follows that if ${\bf a}_1,\ldots ,{\bf a}_{k+1}$ are $k+1 $ $n$-vectors such that
$\bigwedge_{1\leq i\leq k+1}({ f}({\bf a}_i)+1)=1$, then
$\bigvee _{1\leq j<l\leq k+1}({ f}({\bf a}_j{\bf a}_l)+1)=1$ and the proof of Theorem 1 is complete.
\end{proof}

An equivalent form of Claim \ref{claim2} in the proof of Theorem \ref{characterization} is the following lemma of independent interest
\begin{lemma}If $k+1$ subsets $A_i$, $1\leq i\leq k+1$ of a $k$-element set $A$ have odd size,
 then there are $1\leq i,j\leq k+1$, $i\not=j$, such that $A_i\cap A_j$ has odd size.
\end{lemma}
\begin{remark} The number of such pairs can be even. For an example, let k=4, $A:= \{0,1,2,3\}$ and  $A_1,\dots, A_5$ whose corresponding vectors  are  $a_1:=1110$, $a_2:=1101$, $a_3:=0111$, $a_4=1000$,
$a_5=0001$. There are only four  odd intersections, namely  $A_1\cap A_4$, $A_2\cap A_4$, $A_2\cap A_5$ and  $A_3\cap A_5$.
\end{remark}

\end{document}